\documentclass[12pt]{amsart}
\usepackage{graphics}
\usepackage{epsfig}
\usepackage{amsmath}
\usepackage{amsfonts,latexsym}
\usepackage{color}
\usepackage{a4}
\usepackage{url}
\newtheorem{lem}{Lemma}

\newtheorem{prop}{Proposition}

\newcommand{\RR}{\ensuremath{\mathbb{R}}}

\newcommand{\ZZ}{\ensuremath{\mathbb{Z}}}

\DeclareMathOperator{\GL}{GL}
\DeclareMathOperator{\HYP}{HYP}
\DeclareMathOperator{\CUT}{CUT}
\DeclareMathOperator{\Stab}{Stab}
\DeclareMathOperator{\Aut}{Aut}
\DeclareMathOperator{\Sym}{Sym}
\DeclareMathOperator{\Dual}{Dual}
\DeclareMathOperator{\vertt}{vert}

\begin{document}

\author{Mathieu Dutour Sikiri\'c}
\address{Mathieu Dutour Sikiri\'c, Rudjer Boskovi\'c Institute, Bijenicka 54, 10000 Zagreb, Croatia}
\email{mdsikir@irb.hr}

\author{Viatcheslav Grishukhin}
\address{Viatcheslav Grishukhin, CEMI Russian Academy of Sciences, Nakhimovskii prosp.47 117418 Moscow, Russia}
\email{grishuhn@cemi.rssi.ru}

\thanks{Both authors are grateful to the hospitality of the graduate school of mathematics of Nagoya University. The first author is supported by the Croatian Ministry of Science, Education and Sport under contract 098-0982705-2707. The second author is supported by grant of RFBR-CNRS No 05-01-02805}


\title{The decomposition of the hypermetric cone into $L$-domains}

\date{}

\maketitle

\begin{abstract}
The hypermetric cone $\HYP_{n+1}$ is the parameter space
of basic Delaunay polytopes of $n$-dimensional lattice.
If one fixes one Delaunay polytope of the lattice then there are
only a finite number of possibilities for the full
Delaunay tessellations.
So, the cone $\HYP_{n+1}$ is the union of a finite set of $L$-domains,
i.e. of parameter space of full Delaunay tessellations.

In this paper, we study this partition of the hypermetric cone into
$L$-domains.
In particular, we prove that the cone $\HYP_{n+1}$ of
hypermetrics on $n+1$ points contains exactly $\frac{1}{2}n!$ 
principal $L$-domains.
We give a detailed description of the decomposition of
$\HYP_{n+1}$ for $n=2,3,4$ and a computer result for $n=5$
(see Table \ref{TableDataHYPn}).
Remarkable properties of the root system
$\mathsf{D}_4$ are key for the decomposition of $\HYP_5$.
\end{abstract}

\section{Introduction}
\begin{figure}\label{DelaunayVoronoiPolytopes}
\begin{center}
\begin{minipage}{5cm}
\centering
\epsfig{file=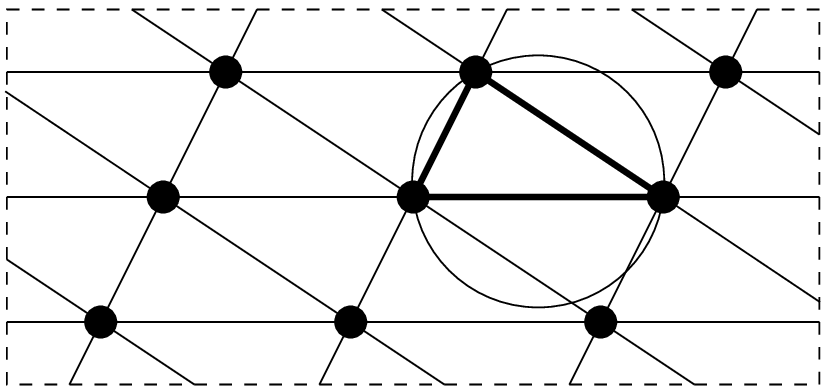,height=20mm}\par
The Delaunay polytopes of $L$
\end{minipage}
\begin{minipage}{5cm}
\centering
\epsfig{file=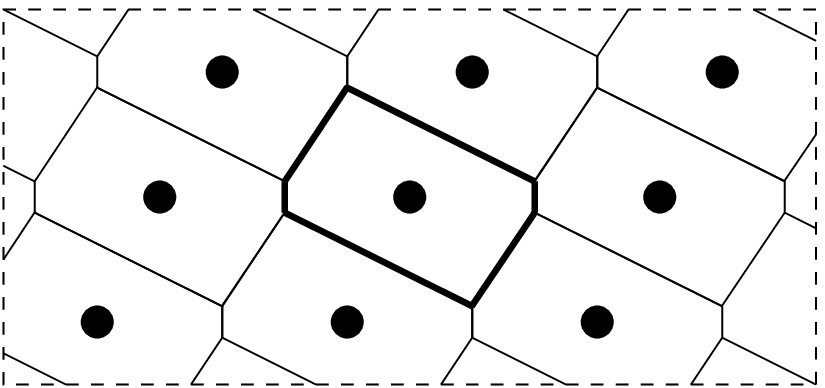,height=20mm}\par
The Voronoi polytope of $L$ and its translates
\end{minipage}
\end{center}
\caption{A lattice $L\subset \RR^2$ and the induced partitions}
\end{figure}

An {\em $n$-dimensional lattice} $L$ is a subgroup of $\RR^n$ of the form
$L=v_1 \ZZ+\dots+v_n \ZZ$ with $(v_1, \dots, v_n)$ a basis of $\RR^n$.
Let $S(c,r)$ be a sphere in $\RR^n$ with center $c$ and radius $r$.
Then, $S(c,r)$ is said to be a {\em Delaunay sphere} in the lattice $L$ if the
following two conditions hold:
\begin{itemize}
\item[(i)] $\Vert v-c\Vert\geq r$ for all $v\in L$,
\item[(ii)] the set $S(c,r) \cap L$ has affine rank $n+1$.
\end{itemize}
The $n$-dimensional polytope $P$, which is defined as the convex hull of
the set $S(c,r)\cap L$, is called a {\em Delaunay polytope} of rank $n$.
The Delaunay polytopes of rank $n$ form a face-to-face tiling of $\RR^n$.
The {\em Voronoi polytope} $P_V(L)$ of a lattice $L$ is the
set of points, whose closest element in $L$ is $0$.
Its vertices are centers of Delaunay polytopes of $L$.
The polytope $P_V(L)$ forms a tiling of $\RR^n$
under translation by $L$, i.e. it is a {\em parallelohedron}
(see Figure \ref{DelaunayVoronoiPolytopes}).

The cones ${\mathcal S}^n_{>0}$, ${\mathcal S}^n_{\geq 0}$ are respectively the
cone of positive definite, positive semidefinite $n\times n$
matrices.
The {\em rational} closure ${\mathcal S}^n_{rat \geq 0}$ of ${\mathcal S}^n_{>0}$
is defined as the positive semidefinite matrices, whose kernel is defined
by rational equalities (see \cite{DSV}).
Given a basis $B=(v_i)_{1\leq i\leq n}$ of a lattice $L$,
we associate the Gram matrix $a=B^T B\in {\mathcal S}^n_{>0}$.
On the other hand if $a\in {\mathcal S}^n_{>0}$, then we can find an
invertible real matrix $B$ such that $a=B^T B$, i.e. $B$ is the basis
of a lattice $L=B \ZZ^n$ with Gram matrix $a$.
So, we can replace the study of Delaunay polytopes of $L$
for the standard scalar product by the study of Delaunay polytopes
of $\ZZ^n$ for the scalar product $x^T a x$.
If one takes another basis $B'$ of $L$, then $B'=BP$ for
some $P\in \GL_n(\ZZ)$ and one has $a'=P^{T} a P$, i.e. $a$ and $a'$
are {\em arithmetically equivalent}.
In other words, the study of $n$-dimensional lattices up to isometric
equivalence is the same as the study of positive definite $n\times n$
symmetric matrices, up to arithmetic equivalence.
In \cite{DSV} it is proved that if $a\in {\mathcal S}^n_{\geq 0}$, then
one can define, possibly infinite, Delaunay polytopes of $\ZZ^n$ for 
$x^T a x$ if and only if $a\in {\mathcal S}^n_{rat \geq 0}$.

Given a polytope $P$ of $\ZZ^n$ the condition that it is
a Delaunay polytope for the norm $x^T a x$ translates to
linear equalities and strict inequalities on the coefficients
of $a$.
An $L$-domain is the convex cone of all matrices 
$a\in {\mathcal S}^n_{>0}$ such that $\ZZ^n$ has the
same Delaunay tessellation for $x^T a x$
(Details see, for example, in \cite{Vo,DL,DSV}).
Voronoi proved that the cone ${\mathcal S}^n_{>0}$
is partitioned into polyhedral $L$-domains.
An $L$-domain of maximal dimension $\frac{1}{2}n(n+1)$ is called {\em primitive}.
An $L$-domain is primitive if and only if the Delaunay tiling related to it consists only of simplices.
Each non-primitive $L$-domain is an open face of the
closure of a primitive one.
In particular, an extreme ray of the closure of an
$L$-domain is a non-primitive one-dimensional $L$-domain.
The group $\GL_n(\ZZ)$ acts on the $L$-domains of ${\mathcal S}^n_{>0}$
by ${\mathcal D}\mapsto P^T {\mathcal D}P$,
and there is a finite number of orbits of $L$-domains, called {\em $L$-types}.
The geometric viewpoint is most useful for thinking,
and drawings about lattice and the Gram matrix viewpoint
is the most suitable to machine computations.

A {\em metric} on the set $\{0,1,\dots,n\}$ is a function $d$
such that $d(x,x)=0$, $d(x,y)=d(y,x)$ and
$d(x,y)\leq d(x,z)+d(z,y)$.
A metric $d$ is an {\em hypermetric} if it satisfies
the inequalities
\begin{equation}\label{hyp}
H_z(d)=\sum_{0\leq i<j\leq n}z_iz_jd(i,j)\leq 0
\end{equation}
for all integral vectors $z\in\ZZ^{n+1}$ such that
$\sum_{i=0}^{n} z_i=1$. The set of all hypermetrics on $n$ points
$\{0,\dots, n-1\}$ is denoted by $\HYP_n$.

The group $\Sym(n)$ acts on $\HYP_n$;
it is proved in \cite{DezaPasechnik} that there is no
other symmetries if $n\not= 4$.
It is proved in \cite{DL} that
$\HYP_{n+1}$ is polyhedral, i.e. among the infinite set of inequalities
of the form (\ref{hyp}), a finite number suffices to get all facets.
This result can be proved in many different ways, see \cite[Theorem~14.2.1]{DL};
the second proof uses that the image $\xi(\HYP_{n+1})$ is the union of
a finite number of $L$-domains.
The purpose of this article is to investigate
such decompositions of $\HYP_{n+1}$.

\begin{table}\label{FacetsOfHYPn}
\begin{center}
\begin{tabular}{|c|c|}
\hline
$n$  & representative of orbits of facets of $\HYP_n$\\
\hline
$3$   & $(1,1,-1)$ (triangle inequality)\\
$4$   & $(1,1,-1,0)$\\
$5$   & $(1,1,-1,0,0)$ and $(1,1,1,-1,-1)$ (pentagonal inequality)\\
$6$   & $(1,1,-1,0,0,0)$, $(1,1,1,-1,-1,0)$,\\
      & $(1,1,1,1,-1,-2)$ and $(2,1,1,-1,-1,-1)$\\
\hline
\end{tabular}
\end{center}
\caption{The facets of $\HYP_n$ for $n\leq 6$}
\end{table}

The set of orbits of facets of $\HYP_{n}$ for $n\leq 6$
is given in Table \ref{FacetsOfHYPn}, 
$\HYP_7$ has $14$ orbits of facets (see \cite{DL,Baranovski,BR})
and the list is not known for $n\geq 8$.
An inequality of (\ref{hyp}) is called {\em $k$-gonal}
if $\sum_{i=0}^{n} |z_i|=k$.
$3$-gonal and $5$-gonal inequalities are also called
{\em triangle} and {\em pentagonal} inequalities, respectively.

A Delaunay polytope $P$ of a lattice $L$ is called {\em generating} if 
the smallest, for the inclusion relation, lattice containing
$V(P)$ is $L$.
Moreover, if there exist a family $(v_0, \dots, v_n)$ of vertices of
$P$ such that
for any $v\in L$ there exist $\alpha_i\in \ZZ$ with
\begin{equation*}
1=\sum_{i=0}^{n} \alpha_i, v=\sum_{i=0}^{n} \alpha_i v_i
\end{equation*}
then $P$ is called {\em basic} and $(v_0, \dots, v_n)$ is an {\em affine
basis}. Given such an affine basis, we define the
distance $d(i,j)=\Vert v_i-v_j\Vert^2$ and we have
\begin{equation*}
H_{b}(d)=\sum_{0\leq i< j\leq n} b_i b_j d(i,j)=(r^2-\Vert \sum_{i=0}^{n} b_i v_i -c\Vert^2)\leq 0, 
\end{equation*}
where $r$ and $c$ relate to the circumscribing sphere $S(c,r)$
of $P$.
So the hypermetric inequalities correspond to the inequalities determining
a family $(v_0, \dots, v_n)$ to be an affine basis of a Delaunay polytope.
Moreover we have $H_{b}(d)=0$ if and only if $\sum_{i=0}^n b_i v_i$ is a
vertex of $P$. In other words the hypermetric cone $\HYP_{n+1}$ is the
parameter space of a basic simplex in $\ZZ^n$.
We refer for proofs to \cite{AS,DL}.
In practice, if $(v_0, \dots, v_n)$ is an affine basis, we can replace it
by $v_0=0$, $v_i=e_i$ and call the corresponding simplex {\em main}.
At this point we should note that the hypermetric cone is just one
possibility for a parameter space of Delaunay polytopes.
Following \cite{E75,E92}, define $M_{2,n}$ to be the space of real
polynomials of $n$ variables with degree at most $2$.
We then have
\begin{equation*}
C_n=\{f\in M_{2,n}\quad |\quad f(x)\geq 0\quad {\rm for}\quad {\rm all}\quad x\in\ZZ^n\}
\end{equation*}
If $P$ is a Delaunay polytope of rank $k\leq n$, then we define
\begin{equation*}
C_n(P)=\{f\in C_n\quad |\quad f(x)=0\quad {\rm for}\quad {\rm all}
\quad x\in \vertt P\}
\end{equation*}
Those cones were used in \cite{newalgo} to find some so-called
{\em perfect} Delaunay polytopes.
Note that if $P=\{0,e_1, \dots, e_n\}$, then the cone $C_n(P)$
is isomorphic to the cone $\HYP_{n+1}$.

The covariance map $\xi:d\to a$ transforms a hypermetric $d$ on $n+1$ points
$i$, $0\le i\le n$, of a set $X$ into an $n\times n$ positive semidefinite
symmetric matrix $a$ as follows:
\[a_{ij}=\xi(d(i,j))=\frac{1}{2}(d(0,i)+d(0,j)-d(i,j)) \]
(see \cite[Section 5.2]{DL}).
The covariance $\xi$ maps the hypermetric cone $\HYP_{n+1}$ into 
${\mathcal S}^n_{rat \geq 0}$.
Note that there are $n+1$ distinct such maps
depending on which point of $\{0,e_1, \dots, e_n\}$
is chosen as the zero point.

\section{Decomposition methods}

Recall that the Delaunay tiling related to a Gram matrix $a$ from a primitive
$L$-domain ${\mathcal D}$ consists of simplices. The set of Delaunay simplices of the
tiling containing the common lattice point $0$ is the star $St_0$. By
translations, along $\ZZ^n$,
the star $St_0$ determines fully the Delaunay tiling of $\ZZ^n$.
The primitive $L$-domain ${\mathcal D}$ belongs to
$\xi(\HYP_{n+1})$ if and only if its star $St_0$ contains a main simplex.

A {\em wall} $W$ is an $\frac{n(n+1)}{2}-1$ dimensional $L$-domain,
which necessarily separates two primitive $L$-domains ${\mathcal D}$,
${\mathcal D}'$.
Let one moves a point $a$ from the primitive $L$-domain ${\mathcal D}$
to ${\mathcal D}'$ by passing through $W$.
When $a\in W$, some pairs of simplices of $St_0$,
which are mutually adjacent by a facet, glue into repartitioning polytopes.
It is well known (see, for example, \cite{DSV,Vo}) that an $n$-dimensional polytope
with $n+2$ vertices can be triangulated in exactly two ways.
When the point $a$ goes from $W$ into
the $L$-domain ${\mathcal D}'$, each repartitioning polytope repartitions into
its other set of simplices.

Since each repartitioning polytope has $n+2$ vertices, there is an affine
dependence between its vertices. This affine dependence generates a linear
equality between the coefficients $a_{ij}$ of the Gram
matrix $a\in {\mathcal S}^n_{>0}$.
This equality is just the equation determining the hyperplane supporting the
wall $W$. If the point $a$ lies inside the $L$-domain $\mathcal D$, then this
equality holds as an inequality.

So, the convex hull of vertices of any pair of simplices of $St_0$ adjacent
by a facet is a putative repartitioning polytope giving an inequality
separating the $L$-domain $\mathcal D$ from another $L$-domain. All adjacent pairs
of simplices of $St_0$ determine a system of inequalities describing the
polyhedral cone of the primitive $L$-domain $\mathcal D$. Note that some of
these inequalities define faces of $\mathcal D$ but not walls. If the adjacent pair of
simplices contains the main simplex then the corresponding wall lies on a
facet of the cone $\xi(\HYP_{n+1})$.

Using the above system of inequalities, one can define all extreme rays of
$\mathcal D$, and then all facets of $\mathcal D$.
The sum of Gram matrices lying on extreme rays of $\mathcal D$
is an interior point $a({\mathcal D})$ of $\mathcal D$ uniquely
related to this $L$-domain (see \cite{DSV} for more details).
Hence, $L$-domains $\mathcal D$ and ${\mathcal D}'$ belong to the
same $L$-type if and only if $a({\mathcal D}')=P^T a({\mathcal D}) P$
for some $P\in \GL_n(\ZZ)$, i.e. $a({\mathcal D}')$ and $a({\mathcal D})$
are arithmetically equivalent.

The algorithm for enumerating primitive $L$-domains in $\xi(\HYP_{n+1})$ works as
follows. One takes a primitive Gram matrix $a\in\xi(\HYP_{n+1})$. There is a
standard algorithm which, for a given $a\in {\mathcal S}^n_{>0}$, constructs its
simplicial Delaunay tiling. (For example, one can take $a$ from a principal
$L$-domain, described in following sections).

Using the star $St_0$  of the Delaunay tiling, the algorithm, for each pair
of adjacent simplices, determines the corresponding inequality.
By the system of obtained inequalities, the algorithm finds all
extreme rays of the domain $\mathcal D$ of $a$, computes the interior central
ray $a({\mathcal D})$ and finds all facets of $\mathcal D$.
The $L$-domain $\mathcal D$ is put in the list ${\mathcal L}$
of primitive $L$-domains in $\xi(\HYP_{n+1})$.

Let $F$ be a facet of $\mathcal D$ which does not lie on a facet of
$\xi(\HYP_{n+1})$.
For each repartitioning polytope related to the facet $F$ of
$\mathcal D$, the algorithm finds another partition into simplices.
This gives the Delaunay tiling of the primitive $L$-domain
${\mathcal D}'$, which is neighboring to $\mathcal D$ by the
facet $F$. The algorithm finds all extreme rays of ${\mathcal D}'$, the
ray $a({\mathcal D'})$ and tests if it is arithmetically equivalent
to $a({\mathcal D})$ for some $\mathcal D$ from ${\mathcal L}$.
If not, one puts ${\mathcal D}'$ in
${\mathcal L}$.
The algorithm stops when all neighboring $L$-domains are equivalent to
ones in ${\mathcal L}$.
This algorithm is very similar to the one in \cite{sturmfels} for the
decomposition of the metric cone into $T$-domains and belongs to the
class of {\em graph traversal algorithms}.

We give some details of the partition of $\HYP_{n+1}$ into
primitive $L$-domains in Table \ref{TableDataHYPn}.

\begin{table}\label{TableDataHYPn}
\begin{center}
\begin{tabular}{|c|c|c|c|c|} \hline
n & \# primitive

&\# facets &\# orbits of primitive & \# primitive $L$-domains   \\
  & $L$-types & of $\HYP_{n+1}$ &$L$-domains in $\HYP_{n+1}$ & in $\HYP_{n+1}$\\
  &    &                 & under $\Sym(n+1)$ & \\ \hline
2 & 1     & 3    & 1     & 1 \\
3 & 1     & 12   & 1     & 3 \\
4 & 3     & 40   & 5     & 172 \\
5 & 222   & 210  & 8287  & 5 338 650 \\
\hline
\end{tabular}
\end{center}
\caption{Decomposition of $\HYP_{n+1}$ into $L$-domains}
\end{table}

We now expose another enumeration method of the orbits of primitive $L$-domains in $\HYP_{n+1}$.
Consider a primitive $L$-domain ${\mathcal D}$.
The group $\Stab({\mathcal D})=\{P\in \GL_n(\ZZ)\mbox{~}:\mbox{~}P^T a({\mathcal D})P=a({\mathcal D})\}$ is a finite group,
which permutes the translation classes
of simplices of the Delaunay decomposition of ${\mathcal D}$.
It splits the translation classes of simplices into different orbits.
Let $S$ be a basic simplex in ${\mathcal D}$; if one choose the coordinates
such that $\vertt S=\{0, e_1, \dots, e_n\}$ then one obtains
an $L$-domain ${\mathcal D}_S$, whose image by $\xi^{-1}$ is
included in $\HYP_{n+1}$.
A permutation of the vertex set of $S$ induces a permutation in $\HYP_{n+1}$
as well.
So, two cones $\xi^{-1}({\mathcal D}_S)$ and $\xi^{-1}({\mathcal D}_{S'})$
are equivalent under $\Sym(n+1)$
if and only if $S$ and $S'$ belong to the same orbit of translation
classes of simplices under $\Stab({\mathcal D})$.
Therefore from the list of $L$-types in dimension $n$, one obtains
the orbits of $L$-domains in $\HYP_{n+1}$.

\section{Dicings, rank $1$ extreme rays of an $L$-domain and of $\HYP_{n+1}$}

We denote by $b^Tc$ the scalar product of column vectors $b$ and $c$.
A vector $v\in \ZZ^n$ is called {\em primitive} if the greatest common
divisor of its coefficients is $1$; such a vector defines a family of
parallel hyperplanes $v^Tx=\alpha$ for $\alpha\in \ZZ$.
In the same way, a {\em vector family} ${\mathcal V}=(v_i)_{1\leq i\leq M}$ 
of primitive vectors defines $M$ families of parallel hyperplanes.
A vector family ${\mathcal V}$ is called a {\em lattice dicing} if for any
$n$ independent vectors $v_{i_1}, \dots, v_{i_n}\in {\mathcal V}$
the vertices of the hyperplane arrangement $v_{i_j}^T x=\alpha_i$
form the lattice $\ZZ^n$
(see an example on Figure \ref{DelaunayVoronoiPolytopes}).
This is equivalent to say that any $n$ independent
vectors $v_i$ have determinant $\pm 1$, i.e. the vector family 
is {\em unimodular}.
Given a dicing, the connected components of the complement of $\RR^n$
by the hyperplane arrangement form a partition of $\RR^n$ by polytopes.
It is proved in \cite{ER} that the polytopes of a lattice dicing
defined by ${\mathcal V}=(v_i)_{1\leq i\leq M}$ are Delaunay polytopes for
the matrices belonging to the $L$-domain generated by the rank
$1$-forms $(v_iv_i^T)_{1\leq i\leq M}$, whose corresponding quadratic
form is $f_{v_i}(x)=(v_i^T x)^2$.
The reverse is also proved there, i.e. any $L$-domain, whose
extreme rays have rank $1$ has its Delaunay tessellation 
being a lattice dicing.
Such $L$-domains are called {\em dicing domains};
they are simplicial, i.e. 
their dimension is equal to their number of extreme rays.

Our exposition of matroid theory is limited here to what is useful for the comprehension
of the paper and we refer to \cite{DaG,Aig,truemper} for more details.
Given a graph $G$ of vertex-set $\{1, \dots, n\}$ we associate to
every edge $e=(i,j)$ a vector $v_e$, which is equal to $1$ in position $i$,
$-1$ in position $j$ and $0$ otherwise.
The vector family ${\mathcal V}(G)=(v_e)_{e\in E(G)}$ is unimodular and
is called the {\em graphic unimodular system of the graph $G$}.
Given a unimodular $n$-dimensional system $U$ of $m$ vectors,
for any basis $B\subseteq U$, we can write $U=B(I_n, A)$, where $A$ is a totally unimodular matrix and $(I_n, A)$ is the concatenation of $I_n$ and $A$.
The matrix $(-A^T, I_{m-n})$ defines a unimodular system,
which is called the {\em dual of $U$} and denoted $\Dual(U)$.
Given a graph $G$ of vertex set $V(G)$ and edge set $E(G)$,
we choose an orientation on every edge $e$ and associate to
it a coordinate $x_{e}$.
We define a vector space $V$ to be the set of vector $v\in \RR^{E(G)}$
satisfying for all $x\in V(G)$ to the {\em vertex cut equation}
\begin{equation*}
0=\sum_{y\in N(x)} v_{(x,y)} \epsilon_{(x,y)}
\end{equation*}
with $N(x)$ the neighbors of $x$ and $\epsilon_{(x,y)}=1$ if the orientation of the edge $(x,y)$ goes from $x$ to $y$ and $-1$ otherwise.
Take $v_1, \dots, v_N$ a basis of the space and denote by $CoGr(G)$ the {\em cographic unimodular system of the graph $G$} defined to be the vector system obtained by taking the transpose of the matrix $(v_1, \dots, v_N)$.
The unimodular systems $CoGr(G)$ and $\Dual(Gr(G))$ are isomorphic.
In \cite{DaG} a general method for describing unimodular vector
families is given using graphic, cographic unimodular systems and
a special unimodular system named $E_5$ (or $R_{10}$ as in \cite{Se}).

Given a finite set $E$, a {\em matroid} $M=M(E)$ is a
family ${\mathcal C}(M)$ of subset of $X$ called
{\em circuits} such that:
\begin{itemize}
\item for $C_1,C_2\in{\mathcal C}(M)$, it holds
$C_1\not\subseteq C_2$, $C_2\not\subseteq C_1$ if $C_1\not= C_2$;
\item if $e\in C_1\cap C_2$, then there is $C_3\in{\mathcal C}(M)$ such that $C_3\subseteq C_1\cup C_2-\{e\}$.
\end{itemize}
A set of vectors $q_e$, $e\in E$, represents a matroid $M(E)$ if, for any
circuit $C\in{\mathcal C}(M)$, the equality $\sum_{e\in C}q_e=0$ holds.
A matroid, is called {\em regular} if it admits a representation as a
unimodular system of vectors.
If $M(E)$ is a {\em graphic} or {\em cographic} matroid of a graph $G$
with a set $E$ of edges, then circuits of $M$ are {\em cycles}
or {\em cuts} of $G$, respectively. 

A graph $G$ is {\em plane} if it is embedded in the $2$-plane
such that any two edges are non-crossings.
A plane graph defines a partition of the plane into faces delimited by
edges. The {\em dual graph} $G^*$ is the graph defined by faces with
an edge between two faces if they share an edge.
Then $(G^*)^*=G$, and there is a
bijection between (intersecting) edges of $G$ and $G^*$ such that each cut
of $G$ corresponds to a cycle of $G^*$, and vice versa. In other words,
the cographic matroid of $G$ and the graphic matroid of $G^*$ are isomorphic.

The only rank $1$ extreme rays of the cone $\HYP_{n+1}$ are {\em cut metrics}.
For $n\le 5$, the hypermetric cone $\HYP_{n+1}$ coincides with the cut cone
$\CUT_{n+1}$ which is the cone hull of $2^{n}-1$ cut metrics. 
Denote $N=\{1,\dots, n\}$; if $S\subset N$, $S\not= \emptyset$
then the cut metric $\delta_S$ on $X=\{0\}\cup N$ is defined as follows:
\begin{equation*}
\delta_S(i,j)=1 \mbox{ if } |\{ij\}\cap S|=1,
\mbox{ and }\delta_S(i,j)=0, \mbox{ otherwise}.
\end{equation*}
The covariance map $\xi$ transforms the cut metric $\delta_S$ into the
following {\em correlation} matrix $p(S)$ of rank 1:
\[p_{ij}(S)=1 \mbox{ if } \{ij\}\subset S, \mbox{ and }
  p_{ij}(S)=0, \mbox{ otherwise, where  }1\le i, j\le n. \]

The quadratic form corresponding to the correlation matrix $p(S)$ is
\[f_S(x)=\xi(\delta_S)(x)=\sum p_{ij}(S)x_ix_j=(\sum_{i\in S}x_i)^2. \]
So, $f_S=f_q$ with $f_q(x)=(q^T x)^2$ and $q=\sum_{i\in S}b_i:=b(S)$
the incidence vector of the set $S$.
In summary:
\begin{lem}\label{bS}
A vector $q\in \ZZ^n$ determines an extreme ray $f_q$ of $\xi(\HYP_{n+1})$ if
and only if $q=b(S)$ for some $S\subseteq N$, $S\not=\emptyset$.
\end{lem}

By Lemma~\ref{bS}, if $U$ determines a dicing domain in $\xi(\HYP_n)$, then
the set of coordinates of vectors from $U$ in the basis $\mathcal B$ forms a
unimodular matrix with $(0,1)$-coefficients.
Note that the columns of any $(0,1)$-matrix are incidence vectors
of subsets of a set.
Since any unimodular system $U$ determines a dicing $L$-domain
${\mathcal D}(U)$, we have the following proposition:
\begin{prop}\label{U01}
Let ${\mathcal D}(U)$ be a dicing domain determined by a unimodular set $U$.
The following assertions are equivalent:

(i) ${\mathcal D}(U)$ lies in $\xi(\HYP_n)$;

(ii) $U$ is represented by a $(0,1)$-matrix.
\end{prop}

\section{The principal $L$-domain}

There is a unique primitive $L$-type, whose $L$-domains are dicing domains of
maximal dimension $\frac{1}{2}n(n+1)$ (\cite{zolotarev77,Dick}).
Voronoi calls this $L$-type
{\em principal}. Each principal $L$-domain is simplicial and all its
$\frac{1}{2}n(n+1)$ extreme rays have rank 1. The set $Q$ of vectors $q$
determining extreme rays $f_q$ of a principal $L$-domain forms a maximal
unimodular system. This system is the classical unimodular root system
$\mathsf{A}_n$ representing the graphic matroid of the complete graph
$\mathsf{K}_{n+1}$ on $n+1$ vertices.

In our case, when a principal $L$-domain is contained in $\xi(\HYP_{n+1})$, its
extreme rays belong to the set $\{p(S):S\subseteq N, S\not=\emptyset\}$ of
extreme rays of the cone $\xi(\HYP_{n+1})$. Hence, the vectors $q$ have
the form $b(S)$ for $S\subseteq N$. We shall find all subsets of these
vectors representing the graphic matroid of $\mathsf{K}_{n+1}$.
We orient edges of $\mathsf{K}_{n+1}$ into {\em arcs} and relate a vector $b(S)$ to each arc of the
directed graph $\mathsf{K}_{n+1}$ such that, for any directed circuit $C$ in
$\mathsf{K}_{n+1}$, the following equality holds
\begin{equation}\label{crc}
\sum_{e\in C}\varepsilon_eb(S_e)=0.
\end{equation}
Here $b(S_e)$ is the vector related to the arc $e$ and $\varepsilon_e=1$ if
the directions of $e$ and $C$ coincide, and $\varepsilon_e=-1$, otherwise.

Given a chain of equally directed arcs labeled by
one-element set, a subchain of this chain determines
a $(0,1)$-characteristic vector.
It is known that the set of characteristic vectors of
a set of connected subchains determines a graphic
unimodular system.
We show below that such graphical systems
are contained in the set $\{b(S):S\subseteq N\}$.

This relation of vectors $b(S)$ and arcs of $\mathsf{K}_{n+1}$ provides a labeling
of arcs of $\mathsf{K}_{n+1}$ by subsets $S\subseteq N$. We call this labeling
{\em feasible} if the corresponding set of vectors $b(S)$ gives a
representation of the graphical matroid of $\mathsf{K}_{n+1}$, i.e. (\ref{crc})
holds for each circuit $C$ of $\mathsf{K}_{n+1}$.

Consider a $k$-circuit $C=\{e_i:1\le i\le k\}$, whose arcs have the same
directions. Suppose that, for $1\le i\le k$, $S_i$ is a label of the arc
$e_i$, and that this labeling is feasible.
Then the equality $\sum_{i=1}^kb(S_i)=0$ holds.
Since the coordinates of the vectors $b(S)$
take $(0,1)$-values, this equality is not possible for a
feasible labeling. Hence, the directed graph $\mathsf{K}_{n+1}$ with a feasible
labeling has no circuit, whose arcs have the same directions.
Any finite directed graph with no circuit has at least one {\em source} vertex
(and a {\em sink} vertex as well).

Now consider a directed $3$-circuit $C=\{e_1,e_2,e_3\}$ of a feasible labeled
$\mathsf{K}_{n+1}$.
Then two arcs of $C$, say the arcs $e_1,e_2$, have directions
coinciding with the direction of $C$, and the third arc $e_3$ has opposite
direction. If $S_i$ is a label of $e_i$,
$i=1,2,3$, then we have the equality $b(S_1)+b(S_2)=b(S_3)$. This equality
is possible only if $S_1\cap S_2=\emptyset$ and $S_1\cup S_2=S_3$. Since any
two adjacent arcs of a complete graph belong to a $3$-circuit, we obtain the
following result:

\begin{lem}
\label{Ce}
Let two arcs $e_i$ and $e_j$ be adjacent in a feasible labeled graph
$\mathsf{K}_{n+1}$ and have labels $S_i$ and $S_j$. Then $S_i\cap S_j=\emptyset$ if
the directions of these arcs coincide in the $2$-path $[e_i,e_j]$. If these
directions are opposite, then either $S_i\subset S_j$ or $S_j\subset S_i$.
\end{lem}

Let $v$ be a source vertex of $\mathsf{K}_{n+1}$ and $E(v)$ be the set of $n$ arcs
incident to $v$. Since all arcs of $E(v)$ go out from $v$, any two arcs
$e,e'\in E(v)$ have opposite directions in their $2$-path $[e,e']$. Let
${\mathcal S}(v)=\{S_e:e\in E(v)\}$. By Lemma~\ref{Ce}, the family ${\mathcal S}(v)$
is a nested family of $n$ mutually embedded distinct subsets. This implies
that the sets $S\in{\mathcal S}(v)$ and the arcs $e\in E(v)$ can be indexed as
$S_i$, $e_i$, $1\le i\le n$, such that $S_i$ is the label of $e_i$ and
$|S_i|=i$.

For $2\le i\le n$, let $g_i$ be the arc of the graph $\mathsf{K}_{n+1}$, which forms
a $3$-circuit with the arcs $e_{i-1},e_i\in E(v)$. Lemma~\ref{Ce} implies that
the arc $g_i$ has the one-element set $S_i-S_{i-1}$ as label, and the
direction of $g_i$ coincides with the direction of $e_{i-1}$ in their $2$-path
$[e_{i-1},g_i]$. Now, it is clear that the $n$ arcs $g_i$ for
$1\le i\le n$, where $g_1=e_1$, form an $n$-path, whose arcs have the
same directions and are labeled by one-element sets.
Recall that a non self-intersecting $n$-path in a graph
with $n+1$ vertices is called a
{\em Hamiltonian path}. We obtain the following result.
\begin{lem}
\label{Ham}
A feasible labeled complete directed graph $\mathsf{K}_{n+1}$ has a Hamiltonian path
such that all its arcs have the same directions and each arc has a
one-element labeling set.
\end{lem}

Let $\{0\}\cup N$ be the set of vertices of $\mathsf{K}_{n+1}$, where the vertex $0$ is
the source. Let $0,i_1,i_2,\dots,i_n$ be the vertices of the Hamiltonian path
$\pi$ in Lemma~\ref{Ham}. The path $\pi$ defines uniquely an orientation and
a feasible labeling of $\mathsf{K}_{n+1}$ as follows. The arc with end-vertices
$i_j,i_k$, where $0\le j < k \le n$, is labeled by the set
$S_{jk}=\{i_{r}:j+1\le r\le k\}\subseteq N$.
If one reverse the above order, then one gets the same family of sets,
and the labeled graph $\mathsf{K}_{n+1}$ gives the same representation
of the unimodular system $\mathsf{A}_n$.
We have

\begin{lem}\label{ord}
Any representation of the graphic matroid of the complete
graph $\mathsf{K}_{n+1}$ by vectors $b(S)$, $S\subseteq N$,
$S\not=\emptyset$, is determined by a complete order of the set $N$.
Two opposite orders determine the same representation.
\end{lem}
Since there are $n!$ complete orders on an $n$-set, as a corollary of
Lemma~\ref{ord}, we obtain our main result.
\begin{prop}\label{An}
The cone $\xi(\HYP_{n+1})$ contains $\frac{1}{2}n!$ distinct principal
$L$-domains.
\end{prop}
So, each principal domain is determined by an order
(and its reverse) of the set $N$.
For the sake of definition, we choose the lexicographically minimal
order $\mathcal O$ from these two orders. Let ${\mathcal S}({\mathcal O})$ be the family
of sets $S\subseteq N$, $S\not=\emptyset$, such that elements of each set $S$
determine a continuous subchain of the $n$-chain, corresponding to the order
$\mathcal O$. A principal domain determined by an order $\mathcal O$ of the set $N$
has $\frac{1}{2}n(n+1)$ extreme rays $p(S)$ for $S\in{\mathcal S}({\mathcal O})$.

Each face $F$ of a dicing $L$-domain in $\xi(\HYP_{n+1})$ is uniquely
determined by its extreme rays $p(S)$, all of rank 1. Set
\[{\mathcal S}(F)=\{S\subseteq N: p(S) \mbox{ is an extreme ray of  }F\}. \]
\begin{prop}
\label{adj}
If $n\ge 4$, then any two principal domains in $\xi(\HYP_{n+1})$ are not
contiguous by a facet.
\end{prop}
\proof If two principal domains
${\mathcal D}({\mathcal O}), {\mathcal D}({\mathcal O}')\subset\xi(\HYP_{n+1})$ share a facet
$F$, then they have $\frac{1}{2}n(n+1)-1$ common extreme rays $p(S)$ for
$S\in{\mathcal S}(F)$. This implies that the families ${\mathcal S}({\mathcal O})$ and
${\mathcal S}({\mathcal O}')$ should differ by one element only.
But, for any two
distinct orders $\mathcal O$ and ${\mathcal O}'$, the families ${\mathcal S}({\mathcal O})$
and ${\mathcal S}({\mathcal O}')$ differ at least by two sets, if $n\ge 4$. For
example, suppose $\mathcal O$ and ${\mathcal O}'$ differ by a transposition of two
elements $i$ and $j$. Then there is at least one subchain in $\mathcal O$
containing $i$ and not containing $j$, which is not a subchain of ${\mathcal O}'$.
The same assertion is true for the order ${\mathcal O}'$. \qed

\section{The decompositions of $\HYP_3$ and $\HYP_4$}
Since, for $n=2$ and $n=3$ there exists only one primitive $L$-type, 
namely, the principal $L$-type,
Proposition~\ref{An} describes completely the decompositions of the cones
$\xi(\HYP_{n+1})$ for $n=2,3$.

Recall that each facet $F$ of $\HYP_{n+1}$ is described by
an inequality (\ref{hyp}).
If $n\le 4$ and $F$ is a facet of $\HYP_{n+1}$, then
$z_i\in\{0,\pm 1\}$. 
Hence, we can denote triangle and pentagonal
facets as $F(ij;k)$ and $F(ijk;lm)$, respectively.
Here $z_i=z_j=1$, $z_k=-1$, $z_l=0$, $l\in X-\{ijk\}$, for 
the triangle facet, and $z_i=z_j=z_k=1$, $z_l=z_m=-1$, $z_r=0$ if
$r\in X-\{ijklm\}$, for the pentagonal facet.

Note that $\xi(\delta_S)=p(S)$.
For $S=\{ij\dots k\}$, set $S=ij\dots k$ and $p(S)=p(ij\dots k)$.

${\bf n=2}$. The cone $\xi(\HYP_3)$ is three-dimensional and simplicial. There
is only one order ${\mathcal O}=(12)$ with
${\mathcal S}({\mathcal O})=\{1,2, 12\}$. Its three extreme rays $p(1), p(2)$ and
$p(12)$ span a principal domain, i.e. $\xi(\HYP_3)$ coincides with a
principal domain.

${\bf n=3}$. The cone $\xi(\HYP_4)$ is six-dimensional and has seven extreme
rays $p(i)$, $i=1,2,3$, $p(ij)$, $ij=12,13,23$, $p(123)$ and 12 facets
$\xi(F(ij;k))$ for $i,j,k\in\{0\}\cup N$. Note that
\[p(1)+p(2)+p(3)+p(123)=p(12)+p(13)+p(23)=p. \]
Hence, the four-dimensional cone
${\mathcal C}_1=\RR_+p(1)+\RR_+p(2)+\RR_+p(3)+\RR_+p(123)$ 
intersects by the ray $\RR_+p$ the three-dimensional cone
${\mathcal C}_2=\RR_+p(12)+\RR_+p(13)+\RR_+p(23)$. The ray
$\RR_+p$ is an interior ray of both.

By Proposition~\ref{An}, $\xi(\HYP_4)$ contains three principal domains.
These three six-dimensional $L$-domains are determined by the three
orders $(123)$, $(132)$ and $(213)$ of $N=\{123\}$.
Denote the domain determined by the order 
$(ijk)$ as ${\mathcal D}_j$, where $j$ is the middle element of the order $(ijk)$. 
Since the three one-element subsets $\{i\}, i\in N$, and the set $N$ give
continuous chains in all the three orders, the four rays $p(1),p(2),p(3)$ and 
$p(123)$ are common rays of all the three domains ${\mathcal D}_1$, ${\mathcal D}_2$ 
and ${\mathcal D}_3$. Hence, the cone ${\mathcal C}_1$ is the common four-dimensional
face of these three principal domains. The domain ${\mathcal D}_i$ is the cone
hull of ${\mathcal C}_1$ and two rays $p(ij)$ and $p(ik)$. The four triangle
facets $\xi(F(jk;i)), \xi(F(jk;0)), \xi(F(0i;j)), \xi(F(0i;k))$ of
$\xi(\HYP_4)$ are also facets of ${\mathcal D}_i$. The other two facets of
${\mathcal D}_i$ separating the domain ${\mathcal D}_i$ from ${\mathcal D}_j$ and
${\mathcal D}_k$ are the cone hulls of ${\mathcal C}_1$ with the rays $p(ij)$ and
$p(ik)$, respectively.

\section{$L$-domains in $\xi(\HYP_5)$}

A {\em parallelohedron} is an $n$-dimensional polytope, whose image under a
translation group forms a tiling of $\RR^n$.
Given a face $F$ of a parallelohedron $P$, the set of faces of $P$ which are translates of $F$ is called the {\em zone} of $P$.
For a parallelohedron $P$ the Minkowski sum $P+z(q)$ may not be a
parallelohedron.
A parallelohedron $P$ is called {\em free along a vector $q$} and the vector
$q$ is called {\em free for a parallelohedron $P$} if the sum $P+z(q)$ is a
parallelohedron (see \cite{Gr1} for more details on this notion).

If $P_q=P+z(q)$ is an $n$-dimensional parallelohedron, then $P_q$
has a {\em non-zero width} along the line $l(q)$ spanned by $q$.
This means that the intersection of $P_q$ with a line parallel to
$l(q)$ is distinct from a point. In this case, the lattice $L_q$ of
the parallelohedron $P_q$ has a {\em lamina} $H$, i.e. a hyperplane
$H$ such that $H$ is transversal to $l(q)$, the intersection
$L_q\cap H$ is an $(n-1)$-dimensional sublattice of $L_q$ and each
Delaunay polytope of $L_q$ lies between two neighbouring layers of
$L_q$ parallel to $L_q\cap H$ (see \cite{closedForms}).
If a Voronoi polytope has a non-zero with along a line $l$, then
the lamina $H$ is orthogonal to $l$.

\subsection{Root lattice $\mathsf{D}_4$}\label{D4section}

The lattice $\mathsf{D}_n$ is defined as
\begin{equation*}
\mathsf{D}_n=\{x\in \ZZ^n\quad | \quad \sum_{i=1}^{n} x_i\equiv 0\pmod 2\}
\end{equation*}
If $\{e_i:1\le i\le n\}$ is an
orthonormal basis of $\ZZ^n$, then the set of shortest vectors of
$\mathsf{D}_n$ is
$\pm e_i\pm e_j$ for $1\le i<j\le n$.
It is the set of all facet vectors of $P_V(\mathsf{D}_n)$
and form an irreducible root system, which we also denote by $\mathsf{D}_{n}$.
There are three translation classes of Delaunay polytopes in $\mathsf{D}_{n}$:
the cross polytope $\beta_n$ whose vertex set is formed by all $e_1\pm e_{i}$
for $1\leq i\leq n$, the half cube $\frac{1}{2} H_n$ whose vertex
set is $\{x\in \{0,1\}^n\quad | \quad \sum_{i=1}^{n} x_i\equiv 0\pmod 2\}$
and a second half cube
$\frac{1}{2}H'_{n}$ whose vertex set is
$\{x\in\{1,2\}\times \{0,1\}^{n-1} \quad |\quad \sum_{i=1}^{n} x_i\equiv 0\pmod 2\}$.
The two half cubes are equivalent under the automorphism group $\Aut(\mathsf{D}_n)$ of the root lattice $\mathsf{D}_n$.
It is proved in \cite{Gr,DG} that the Voronoi polytope
$P_V(\mathsf{D}_n)$ of the $n$-dimensional
root lattice $\mathsf{D}_n$ is free only along vectors which
are parallel to edges of $P_V(\mathsf{D}_n)$. Thus there are $2^n+2n$ free vectors.

It turns out that when $n=4$ all Delaunay polytopes are isometric to the cross-polytope $\beta_4$ and equivalent under $\Aut(\mathsf{D}_4)$.
Any $2$-face of $\beta_4$ is contained in three Delaunay polytopes
$\beta_4$ of the Delaunay tessellation.
This proves that the $L$-domains of $\mathsf{D}_4$
are not dicing domains.
Furthermore, $\mathsf{D}_4$ is a {\em rigid lattice}
(see \cite{rigidlattice}), i.e. its 
Delaunay tessellation determines the Gram matrix up to a scalar multiple.
This means that $\mathsf{D}_4$ determine a $1$-dimensional $L$-type
and that any $L$-domain containing it as an extreme ray is not
a dicing domain.

The free vectors of $P_{V}(\mathsf{D}_4)$ are parallel and can
be identified to the diagonals of the cross polytopes. They are
\begin{itemize}
\item $\pm 2e_i$ with $1\leq i\leq 4$ for $\beta_4$
\item $(\pm 1, \pm 1, \pm 1, \pm 1)$ with even plus signs for $\frac{1}{2}H_4\equiv \beta_4$
\item $(\pm 1, \pm 1, \pm 1, \pm 1)$ with odd plus signs for $\frac{1}{2}H'_4\equiv \beta_4$
\end{itemize}
Up to a factor $\sqrt{2}$, those $24$ vectors are an isometric copy of the root system $\mathsf{D}_4$, which we denote by $\mathsf{D}_{4,2}$. The  union $\mathsf{D}_4 \cup \mathsf{D}_{4,2}$ is the irreducible root system $\mathsf{F}_4$ (see \cite{Humph}).

The Voronoi polytope $P_V(\mathsf{D}_4)$ of the root lattice $\mathsf{D}_4$
is the regular polytope $24$-cell, whose automorphism group is the Coxeter
group $W(\mathsf{F}_4)$ of Schl\"afli symbol $\{3,4,3\}$ (see \cite{Cox}).
Its number of vertices, $2$-faces, $3$-faces is $24$, $96$, $24$.
Each facet is an octahedron with four pairs of opposite and
mutually parallel triangular $2$-faces. The facet vectors of the Voronoi
polytope $P_V(\mathsf{D}_4)$ are the $24$ roots of the root system
$\mathsf{D}_4$.
The polytope $P_V(\mathsf{D}_4)$ has $12$ edge {\em zones} of mutually
parallel edges representing $\mathsf{D}_{4,2}$ and
$16$ face zones of mutually parallel triangular faces.
Each edge zone contains $8$ parallel edges, and each face zone contains six
parallel faces.

We choose a basis of $\mathsf{D}_4$ and denote by $a(\mathsf{D}_4)$ the
Gram matrix of $\mathsf{D}_4$ in this basis.
The $24$ vertices of $P_V(\mathsf{D}_4)$ are given by $\frac{1}{2}\mathsf{D}_{4,2}$.

\subsection{The $L$-domains containing $a(\mathsf{D}_4)$}
For $n=4$, there are three primitive $L$-types of four-dimensional lattices:
the principal type, and $L$-types called by Delaunay in \cite{De} types II and III.
The ten-dimensional $L$-domains of these $L$-types are constructed as follows (cf., \cite{DG}).

Each $k$-dimensional face of a principal
domain relates to the graphic matroid of a subgraph
$G\subseteq \mathsf{K}_5$ on $k$ edges.
Hence each facet (of dimension $9$) of a
principal domain relates to the graphic matroid of 
$\mathsf{K}_5-{\bf 1}$, i.e. the complete graph $\mathsf{K}_5$ without
one edge.
The Gram matrix $a(\mathsf{D}_4)$ of the root lattice $\mathsf{D}_4$
is an extreme ray of $L$-domains of type II and III.
The cone hull of a facet of a principal domain and of a ray of
type $a(\mathsf{D}_4)$ is an $L$-domain of type II.
Hence, any principal domain is contiguous in ${\mathcal S}^4_{>0}$
by facets only with $L$-domains of type II.

An $L$-domain of type II has the following three types of facets.
One dicing facet by which it is contiguous in ${\mathcal S}^4_{>0}$
to a principal domain relates to the graphic matroid of the graph
$\mathsf{K}_5-{\bf 1}$.
Each of two other types of facets is the cone hull of the ray of
type $a(\mathsf{D}_4)$ and a dicing $8$-dimensional face related
to the graphic matroids of $\mathsf{K}_5-2\times{\bf 1}$ or
of $\mathsf{K}_5-{\bf 2}$.
Here each of these graphs is the complete graph $\mathsf{K}_5$
without two non-adjacent or two adjacent edges, respectively.

The complete bipartite graph $\mathsf{K}_{ij}$ is formed by two blocks $S_1$, $S_2$ of vertices with $|S_1|=i$, $|S_2|=j$ and two vertices adjacent if and only if they belong to different blocks.
An $L$-domain of type III is the cone hull of $a(\mathsf{D}_4)$ and
a $9$-dimensional dicing facet related to the cographic matroid
$CoGr(\mathsf{K}_{33})$ of the bipartite graph $\mathsf{K}_{33}$.
Each $8$-element submatroid of $CoGr(\mathsf{K}_{33})$ is graphic and
relates to the graph $\mathsf{K}_5-2\times{\bf 1}$.
Hence, each other facet of an $L$-domain of type III is
the cone hull of $a(\mathsf{D}_4)$ and a dicing $8$-dimensional
face related to $\mathsf{K}_5-2\times{\bf 1}$.
In ${\mathcal S}^4_{>0}$, this facet is a common facet of
$L$-domains of types II and III.

So, an $L$-domain of type II is contiguous in ${\mathcal S}^4_{>0}$
to $L$-domains of all three types.
It is useful to note (\cite{RB,DV}) that if $f$ belongs to the
closure ${\mathcal D}$
of an  $L$-domain of type II or III then the Voronoi polytope
$P_V(f)$ of $\ZZ^n$ under the quadratic form $f$ is an affine image of
the Minkovski sum
\[\sum_{q\in U}\lambda_qz(q)+\lambda P_V(\mathsf{D}_4), \mbox{  }
\lambda_q\ge 0, \lambda\ge 0, \]
where $U$ is the unimodular set of vectors related to rank $1$ extreme
rays of $\mathcal D$, and $P_V(\mathsf{D}_4)$ is the Voronoi polytope
of the root lattice $\mathsf{D}_4$, whose form $a(\mathsf{D}_4)$ lies
also on extreme ray of $\mathcal D$.

\subsection{Unimodular systems in $\mathsf{D}_{4,2}$}
Let $D(4)\subset \mathsf{D}_{4,2}$ be a subset of $12$ roots chosen by one
from each pair of opposite roots.
The vector system $D(4)$ is partitioned into three
disjoint {\em quadruples} $Q_i$, $i=1,2,3$, of mutually orthogonal roots,
given in Subsection \ref{D4section}.

For what follows, we have to consider {\em triples} $(r_1,r_2,r_3)$ of roots
chosen by one from each quadruple $Q_i$, i.e. $r_i\in Q_i$, $i=1,2,3$.
Let $t=(r_1,r_2,r_3)$ be such a triple. Let $\{ijk\}=\{123\}$, i.e. these
three indices are distinct. 
We have $r_i^2=2$, $r_i^Tr_j\in\{\pm 1\}$.
Hence, the vectors
\begin{equation}\label{rij}
r_{ij}=r_i-(r_i^Tr_j)r_j \mbox{ for }ij=12,23,31,
\end{equation}
are roots of $\mathsf{D}_{4,2}$. Since $r^T_{ij}r_i=-r^T_{ij}r_j=1$, one of
two opposite roots $\pm r_{ij}$ belongs to $Q_k$, say $r_{ij}\in Q_k$.
Hence $r_{12}\in Q_3, r_{23}\in Q_1, r_{31}\in Q_2$.
We have two cases:
\begin{itemize}
\item[(i)] $r_{ij}$ belongs, up to sign, to the triple $t$, i.e. $r_{ij}=r_k$,
for all pairs $ij$;
\item[(ii)] $r_{ij}$ does not belong to $t$, i.e. $r_{ij}\not=r_k$.
\end{itemize}
In case (i), the vectors $r_i$, $i=1,2,3$, are linearly dependent, and
the triple $t$ spans a $2$-dimensional plane.
We say that the triple $t$ is of {\em rank} $2$.
Note that any two roots $r,r'$ from distinct quadruples determine uniquely
the third root $r''=r-(r^Tr')r'$ and that there are $16$ distinct triples
of rank $2$.

In case (ii), the roots $r_i$ of the triple $t$ are linearly independent.
We say that the triple $t$ has rank $3$. In this case the roots 
$r_{ij}$ for $ij=12,23,31$, are distinct and do not coincide with the roots 
$r_i$, $i=1,2,3$.
Moreover, it is not difficult to verify that the triple
$(r_{12},r_{23},r_{31})$ has rank $2$.

Triples of rank $2$ and $3$ are realized in $\sqrt{2}P_V(\mathsf{D}_4)$ as follows. The
roots of a triple $t$ of rank $2$ are parallel to edges of a $2$-face of
$\sqrt{2}P_V(\mathsf{D}_4)$.
We say that a triple $t$ of rank $2$ forms a face of
$\sqrt{2}P_V(\mathsf{D}_4)$.
Since a triple of rank $2$ forms a face of $\sqrt{2}P_V(\mathsf{D}_4)$,
the $16$ triples of rank $2$ relate to the $16$ zones of
triangular faces of it.

Let $t'=(r'_1,r'_2,r'_3)$ be a triple of rank $3$. The 6 vectors $\pm r'_i$,
$i=1,2,3$, have end-vertices in vertices of $\sqrt{2}P_V(\mathsf{D}_4)$. From each
pair $\pm r'_i$ of opposite roots, we choose a vector $r_i$ such that
$r_i^Tr_j=1$ for $ij=12,23,31$. Then end-vertices of the roots $r_i$ are
vertices of a face $F$ of $\sqrt{2}P_V(\mathsf{D}_4)$, and the triple
$t=(r_{ij}=r_i-r_j:ij=12,23,31)$ has rank $2$ and forms the face $F$.

For what follows, we need subsets $U$ of $D(4)$ which are maximal by inclusion
such that $P_V(\mathsf{D}_4)+\sum_{q\in U}\lambda_qz(q)$ is a parallelohedron.
Note that each quadruple $Q_i$ is a basis of $\RR^4$.
Obviously, it is a unimodular set.
It is a maximal by inclusion unimodular subset of $D(4)$,
since any other vector of $D(4)$ has half-integer coordinates in this basis.
However, it is proved in \cite{DG}, that $P_V(\mathsf{D}_4)+\sum_{q\in U}\lambda_qz(q)$
is not a parallelohedron if $U\subset D(4)$ is a quadruple, and it is a
parallelohedron for any other maximal unimodular subsets $U\subseteq D(4)$.
Of course, a maximal unimodular set $U\not=Q_i$ does not contain each
quadruple $Q_i$ as a subset.

\begin{figure}
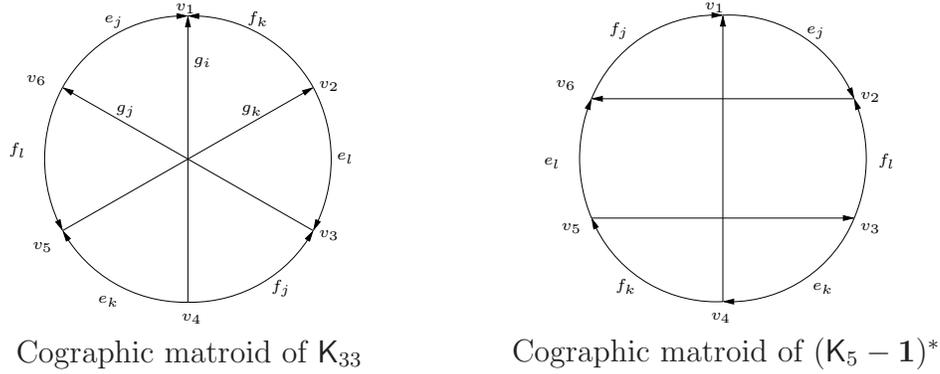

\begin{center}
\begin{minipage}{7cm}
\centering
\input{cographicK33.pstex_t}\par
Cographic matroid of $\mathsf{K}_{33}$
\end{minipage}
\begin{minipage}{7cm}
\centering
\input{cographicDualK5m1.pstex_t}\par
Cographic matroid of $(\mathsf{K}_5-{\bf 1})^*$
\end{minipage}
\end{center}
\caption{Two cographic matroids}
\label{FigureCographicMatroids}
\end{figure}

We show below that maximal unimodular subsets in $D(4)$ represent either the
graphic matroid of the graph $\mathsf{K}_5-{\bf 1}$ or the cographic
matroid of the
graph $\mathsf{K}_{33}$ (see Figure \ref{FigureCographicMatroids}).

The graph $\mathsf{K}_5-{\bf 1}$ is planar. Hence, the graphic matroid of
$\mathsf{K}_5-{\bf 1}$ is isomorphic to the cographic matroid of the graph
$(\mathsf{K}_5-{\bf 1})^*$.
Both graphs $(\mathsf{K}_5-{\bf 1})^*$ and $\mathsf{K}_{33}$ are cubic
graphs on six vertices and nine edges. These graphs have a Hamiltonian
$6$-circuit $C_6$ on six vertices $v_i$, $1\le i\le 6$. Call $C_6$ with the
six edges $(v_i,v_{i+1})$ by a {\em rim}, and the other three edges by
{\em spokes}. The spokes are the following edges: $e_1=(v_1,v_4)$,
$e_2=(v_2,v_6)$, $e_3=(v_3,v_5)$ in $(\mathsf{K}_5-{\bf 1})^*$, and
$e_i=(v_i,v_{i+3})$, $i=1,2,3$, in $\mathsf{K}_{33}$
(see Figure \ref{FigureCographicMatroids}).
Note that the edges $e_2$ and
$e_3$ do not intersect in $(\mathsf{K}_5-{\bf 1})^*$ and intersect in $\mathsf{K}_{33}$. The
described form of $\mathsf{K}_{33}$ is the graph $Q_4$ of \cite{DaG}.

Besides the planarity, the graph $(\mathsf{K}_5-{\bf 1})^*$ differs from $\mathsf{K}_{33}$
by the number of cuts of cardinality three.
All cuts of cardinality $3$ of $\mathsf{K}_{33}$ are the six
one-vertex cuts, i.e. cuts containing three edges
incident to a vertex.
The graph $(\mathsf{K}_5-{\bf 1})^*$, besides the six one-vertex cuts,
has a separating cut containing three non-adjacent edges one of which
is a spoke. (These are the edges $(v_2,v_3)$, $(v_5,v_6)$ and the
spoke $(v_1,v_4)$ of the above description of $(\mathsf{K}_5-{\bf 1})^*$.)

\begin{prop}\label{max}
A maximal by inclusion unimodular subset $U$ of $D(4)$, $U\not=Q_i$,
$i=1,2,3$, is obtained by a deletion from $D(4)$ of a triple $t=(r_1,r_2,r_3)$ where
$r_i\in Q_i$, $i=1,2,3$. Then

(i) if $t$ has rank $2$, then $U$ represents the cographic matroid of the
graph $\mathsf{K}_{33}$;

(ii) if $t$ has rank $3$, then $U$ represents the graphic matroid of the graph
$\mathsf{K}_5-{\bf 1}$, which is isomorphic to the cographic matroid of the dual
graph $(\mathsf{K}_5-{\bf 1})^*$.
\end{prop}
\proof Since $Q_i\not\subseteq U$, for $i=1,2,3$, $U$ does not contain
a triple $t$. We show that a deletion from $D(4)$ of any triple gives a
unimodular set.

Denote by $V(r)$ the set of all triples of rank $2$
containing a root $r\in D(4)$. If $r\in Q_i$ and $r'\in Q_j$, then
$V(r)\cap V(r')$ is the unique triple of rank $2$ containing $r$ and $r'$
if $i\not=j$, and $|V(r)\cap V(r')|=0$ if $i=j$ and
$r\not=r'$.

Let $r_i\in Q_i$, $i=1,2,3$, be the roots of the deleted triple $t$.
Then
\[|V(r_i)|=4, \mbox{  }|V(r_i)\cap V(r_j)|=1, \mbox{ and }
|V(r_1)\cap V(r_2)\cap V(r_3)|=0 \mbox{ or }1. \]
Here, whether $0$ or $1$ stays in the last equality depends on the triple $t$
has rank $3$ or $2$, respectively.
By the inclusion-exclusion principle, we have
\[|\cup_{i=1}^3V(r_i)|=\sum_{i=1}^3|V(r_i)|-\sum_{1\le i<j\le 3}
|V(r_i)\cap V(r_j)|+|\cap_{i=1}^3V(r_i)|. \]
Hence, we have the following two cases.

(i) If $|\cap_{i=1}^3V(r_i)|=1$, i.e. if the triple $t$ has rank $2$, then
$|\cup_{i=1}^3V(r_i)|=3\cdot 4-3\cdot 1+1=10$.

(ii) If $|\cap_{i=1}^3V(r_i)|=0$, i.e. if the triple $t$ has rank $3$, then
$|\cup_{i=1}^3V(r_i)|=9$.

Let $V=\cup_{r\in D(4)}V(r)-\cup_{i=1}^3 V(r_i)$.
Then $V$ is the set of triples of rank 2 that do not
contain the deleted roots $r_1,r_2$ and $r_3$.
Using the set $V$ we construct a cubic directed
graph $G$ as follows.
We take $R_0=D(4)-t$ to be the set of $9$ arcs of $G$.
As a set of vertices of $G$ we take a subset
$V_0\subseteq V$ of six triples.
An arc $r\in R_0$ is incident to a vertex $v\in V_0$ if the root $r$ belongs
to the triple $v$ of rank $2$.
Hence, each vertex $v$ is incident to three edges, and
the graph $G$ is cubic. Directions of arcs are chosen as follows.
Let the roots $r,s,p$ of a triple $v$ satisfy the equation
$r-(r^Ts)s-(r^Tp)p=0$.
Suppose that the
arc $r$ comes in $v$. Then the arc $s$ goes out or comes in the vertex $v$
if $r^Ts=1$ or $r^Ts=-1$, respectively. The same assertion is true for the
arc $p$.

We show that each root $r\in R_0$ belongs to exactly two triples of $V_0$.
Let $r\in Q_k$, then $|V(r)\cap V(r_k)|=0$ and $|V(r)\cap V(r_i)|=
|V(r)\cap V(r_j)|=1$, where $\{ijk\}=\{123\}$. We have two cases:
\begin{equation}\label{neq}
V(r)\cap V(r_i)=V(r)\cap V(r_j) \mbox{ or }
V(r)\cap V(r_i)\not=V(r)\cap V(r_j).
\end{equation}

Suppose that the inequality in (\ref{neq}) holds for $r$.
Since $|V(r)|=4$, only two triples from $V(r)$ belong to $V$.

Now, if the equality in (\ref{neq}) holds, then it implies that
$r=(r^Tr_i)r_i+(r^Tr_j)r_j$. If $t$ has rank $2$, then, up to sign, this
gives $r=r_k$, which contradicts to $r\in R_0$. Hence, if $t$ has
rank $2$, the inequality of (\ref{neq}) holds for all $r\in R_0$. We can set
$V_0=V$, since $V$ contains $16-10=6$ triples.
It is easy to verify that $G$ is isomorphic to $\mathsf{K}_{33}$.

If $t$ has rank $3$, then $r=r_{ij}$, where $r_{ij}$ is defined in equation (\ref{rij}).
The roots $r_{ij}$, $ij=12,23,31$ when the equality holds in \ref{neq},
form a triple $v_0$ of rank $2$.
Obviously, $v_0\in V$. We set $V_0=V-\{v_0\}$. Since
the root $r_{ij}$ belongs to two triples $v_0$ and $(r_{ij},r_i,r_j)$,
the remaining triples of $V(r_{ij})$ belong to $V_0$.
We saw that each root $r\in R_0$, $r\not=r_{ij}$
for $ij=12,23,31$, belongs to two triples of $R_0$.
Hence, in the case when $t$ has rank $3$, the graph $G$ is well defined.
It is easy to verify that $G$ is isomorphic to $(\mathsf{K}_5-{\bf 1})^*$
and $v_0$ corresponds to a separating cut of cardinality $3$.

Arcs of $G$ are labeled naturally by roots from the set $R_0$.
This labeling gives a representation of the cographic matroid
of $G$ by vectors of $R_0$.

Note that triples of rank $2$ are equivalent under action of the automorphism
group of $P_V(\mathsf{D}_4)$. Similarly, all triples of rank $3$ are equivalent under
the automorphism group of $P_V(\mathsf{D}_4)$ extended by changing signs of roots.
Hence, any explicit representations of the cographic matroids of the graphs
$\mathsf{K}_{33}$ and $(\mathsf{K}_5-{\bf 1})^*$ for fixed triples of rank $2$ and $3$ proves that
above labeling give representations for all pairs of triples of rank $2$ and $3$. \qed

\subsection{$L$-types in $\xi(\HYP_5)$}
From Proposition \ref{max} we deduce that each $L$-domain corresponding
to $\mathsf{D}_4$ is contained in $64=4^3$ different $L$-types.
$48$ of them are of type II and $16$ are of type III.

Recall that $P_V(\mathsf{D}_4)$ is free along lines spanned by roots of the
root system $\mathsf{D}_{4,2}$, vectors of which are parallel to diagonals of
the cross-polytopes $\beta_4$ of the Delaunay partition of the lattice
$\mathsf{D}_4$.
Note that bases related to the forms $a(\mathsf{D}_4)\in\xi(HYP_5)$
contain a diagonal of a $\beta_4$.

\begin{lem}\label{lam}
Let the basis related to $a(\mathsf{D}_4)$ contains a diagonal $q\in \mathsf{D}_{4,2}$
of a cross-polytope $\beta_4$. Then the $L$-domain of the parallelohedron
$P_q=P_V(\mathsf{D}_4)+z(q)$ related to this basis, i.e. the $L$-domain of the
form $a(\mathsf{D}_4)+\lambda f_q$, does not belong to $\xi(HYP_5)$.
\end{lem}
\proof The parallelohedron $P_q$ has a non-zero width along the
line $l(q)$ parallel to $q$. The lattice $L_q$ of $P_q$ has a lamina
$H$ which is orthogonal to $q$. The lamina $H$ separate the
cross-polytope $\beta_4$ with a diagonal $q$ into two Delaunay
polytopes, each being a pyramid with a base $\beta_3$ orthogonal to $q$
and lying in $H$. These two pyramids have the end-points of the diagonal
$q$ as apexes. Hence vertices of the basic simplex of
$a(\mathsf{D}_4)+\lambda f_q$ belong to two distinct Delaunay polytopes. This
implies that the $L$-domain of $L_q$ does not belong to $\xi(HYP_5)$. \qed

\begin{prop}\label{main}
The cone $\xi(\HYP_5)$ contains $12$ principal $L$-domains,
$120$ $L$-domains of type II and $40$ $L$-domains of type III,
total $172$ $L$-domains.
\end{prop}
\proof By Proposition~\ref{An}, $\xi(\HYP_{4+1})$ contains
$\frac{1}{2}4!=12$ principal $L$-domains.

The closure of an $L$-domain of types II and III is the convex
hull of $a(\mathsf{D}_4)$ and a dicing facet $F(U)$ related
to a unimodular subset $U\subset \mathsf{D}_{4,2}$
of $9$ vectors that are free for $P_V(\mathsf{D}_4)$.

Each tile of the Delaunay tiling of the root lattice $\mathsf{D}_4$
is the regular four-dimensional cross polytope $\beta_4$.
Any affine base of $\beta_4$ contains exactly the two vertices $w, w'$
of a diagonal and
$3$ vertices $v_1, v_2, v_3$ chosen from the other diagonals.
We have $d_{\beta_4}(w,w')=4$ and $2$ for all other pairs.
There are ${5\choose 2}=10$ ways to choose a pair $\{w,w'\}$ in a $5$
elements set so there are exactly $10$ rays $a_i(\mathsf{D}_4)$
representing $\mathsf{D}_4$ in $\HYP_5$.

Every $L$-domain $a_i(\mathsf{D}_4)$ is contained in the closure
of $64$ primitive $L$-domains ${\mathcal D}(U)$,
but not all of them are included in $\xi(\HYP_5)$.
Each $L$-domain ${\mathcal D}(U)$ has the form
${\mathcal D}(U)=conv(a_i(\mathsf{D}_4)+F(U))$,
where the subset $U\subset \mathsf{D}_{4,2}$ is obtained by a deletion of a
triple from $\mathsf{D}_{4,2}$. By Lemma~\ref{lam}, the inclusion
${\mathcal D}(U)\subset \xi(HYP_5)$ implies $U$ does not contain the
diagonal $q$ which is a basic vector of $a_i(\mathsf{D}_4)$.
Let $q=r_1\in Q_1$, then we have $4$
choices for the root $r_2$ in $Q_2$. For the third root
either we have a rank $2$ triple and $r_3$ is completely determined
or we have a rank $3$ triple in which case there are $3$ choices.
This means that a ray $a_i(\mathsf{D}_4)$ is contained in $4$
$L$-domains of type II included in $\xi(\HYP_5)$
and $12$ $L$-domains of type III included in $\xi(\HYP_5)$.
This means that there are $120$, respectively $40$ $L$-domains of
type II, respectively III in  $\xi(\HYP_5)$. \qed


\begin{thebibliography}{99}

\bibitem[Aig79]{Aig}
M. Aigner, {\em Combinatorial Theory}, Springer-Verlag, Berlin, New-York,
1979.

\bibitem[As82]{AS}
P. Assouad, {\em Sous-espaces de $L^1$ et in\'egalit\'es hyperm\'etriques},
Compte Rendus de l'Acad\'emie des Sciences de Paris,
{\bf 294(A)} (1982) 439--442.

\bibitem[Ba99]{Baranovski}
E.P. Baranovskii, {\em The conditions for a simplex of $6$-dimensional lattice
to be $L$-simplex}, (in Russian) Nauchnyie Trudi Ivanovo state university.
Mathematica, {\bf 2} (1999) 18--24.

\bibitem[BaGr01]{rigidlattice}
E.P. Baranovskii and V.P. Grishukhin,
{\em Non-rigidity degree of a lattice and rigid lattices},
European J. Combin. {\bf 22} (2001) 921--935.




\bibitem[Cox48]{Cox}
H.S.M. Coxeter, {\em Regular polytopes}. Third edition. Dover Publications,
Inc., 1973.

\bibitem[DG99]{DaG}
V.I. Danilov and V.P. Grishukhin,
{\em Maximal unimodular systems of vectors},
European Journal Combinatorics {\bf 20} (1999) 507--526.

\bibitem[De29]{De}
B.N. Delaunay,
{\em Sur la partition r\'eguli\`ere de l'espace \`a 4 dimensions},
Izvestia AN SSSR, ser. phys.-mat. {\bf 1} (1929) 79--110,
{\bf 2} (1929) 145--164.

\bibitem[DGP06]{DezaPasechnik}
A. Deza, B. Goldengorin and D.V. Pasechnik, 
{\em The isometries of the cut, metric and hypermetric cones}, 
Journal of Algebraic Combinatorics {\bf 23-2} (2006) 197--203.

\bibitem[DeG02]{closedForms}
M. Deza and V.P. Grishukhin, 
{\em Rank $1$ forms, closed zones and laminae}, 
Journal de Th\'eorie des Nombres de Bordeaux
{\bf 14} (2002) 103--112.

\bibitem[DL97]{DL}
M. Deza and M. Laurent,
{\em Geometry of Cuts and Metrics},
Springer-Verlag, 1997.

\bibitem[DeG03]{DG}
M. Deza and V.P. Grishukhin,
{\em Once more about 52 four-dimensional parallelotopes},
Taiwanese J. of Math. {\bf 12} (2008) 901--916.

\bibitem[Di72]{Dick}
T.J. Dickson, {\em On Voronoi reduction of positive definite quadratic forms}, 
J. Number Theory {\bf 4} (1972) 330--341.

\bibitem[DR03]{newalgo}
M. Dutour and K. Rybnikov,
{\em A new algorithm in geometry of numbers},
proceedings of the $4^{th}$ international symposium on Voronoi Diagrams
in Science and Engineering, to appear.

\bibitem[DV05]{DV}
M. Dutour and F. Vallentin, {\em Some six-dimensional rigid forms}
in: Voronoi's impact on modern science,
Book 3. Proc. of the third Conf. on Analytic Number
Theory and Spatial Tesselations. Inst. of Math. Kyiv 2005, pp.102--108.

\bibitem[DSV06]{DSV}
M. Dutour, A. Sch\"{u}ermann and F. Vallentin, {\em A generalization of
Voronoi's reduction theory and its application},
Duke Mathematical Journal {\bf 142-1} (2008) 127--164.

\bibitem[ER94]{ER}
R. Erdahl and S. Ryshkov, {\em On lattice dicings},
European Journal of Combinatorics {\bf 15} (1994) 459--481.

\bibitem[Er75]{E75}
R. Erdahl, {\em A convex set of second-order inhomogeneous polynomials
with applications to quantum mechanical many body theory}, Mathematical 
Preprint \#1975-40, Queen's University, Kingston, Ontario.

\bibitem[Er92]{E92}
R. Erdahl, A cone of inhomogeneous second-order polynomials,
{\em Discrete Comput. Geom}. {\bf 8-4} (1992) 387--416.


\bibitem[Gr04]{Gr}
V.P. Grishukhin, {\em Parallelotopes of non-zero width}, Math. Sbornik
{\bf 195-2} (2004) 59--78, (translated in: Sbornik: Mathematics {\bf 195-5}
(2004) 669--686).

\bibitem[Gr06]{Gr1}
V.P. Grishukhin, {\em Free and non-free Voronoi polytopes}, Math. Zametki
{\bf 80-3} (2006) 367--378.

\bibitem[Hu90]{Humph}
J. Humphreys, {\em Reflection groups and Coxeter groups}, Cambridge studies in advanced mathematics {\bf 29}, 1990.


\bibitem[KoZo77]{zolotarev77}
A.N. Korkine and E.I. Zolotarev,
{\em Sur les formes quadratiques positives},
Math. Ann. {\bf 11} (1877) 242--292.



\bibitem[RB]{BR}
S.S. Ryshkov and E.P. Baranovskii, {\em The Repartitioning Complexes in n-dimensional Lattices (with Full Description for $n\leq 6$)}, in Voronoi's impact on modern science, Book 2, Institute of Mathematics, Kyiv (1998) 115-124.


\bibitem[RB05]{RB}
S.S. Ryshkov and E.A. Bolshakova, {\em Theorie of fundamental parallelohedra}
Izvestia RAN, Ser. math. {\bf 69-6} (2005) 187--210.
(Translated in: Izvestia Math. {\bf 69-6} (2006) 1257-1277)



\bibitem[Se80]{Se}
P.D. Seymour, {\em Decomposition of regular matroids}, J. Comb. Theory 
ser. {\bf B} {\bf 28} (1980) 305--359. 


\bibitem[SY04]{sturmfels}
B. Sturmfels and J. Yu, {\em Classification of six-point metrics}, 
Electron. J. Combin. {\bf 11} (2004), Research paper R44.

\bibitem[Tr92]{truemper}
K. Truemper, {\em Matroid decomposition}, Academic Press, 1992.


\bibitem[Vo08]{Vo}
G.F. Voronoi,
{\em Nouvelles applications des param\`etres continus \`a la th\'eorie
de formes quadratiques - Deuxi\`eme m\'emoire.} J. f\"ur die reine und
angewandte Mathematik, {\em 1.Recherches sur les parall\'elo\`edres
primitifs} {\bf 134} (1908) 198--287; {\em 2. Domaines de formes
quadratiques correspondant aux diff\'erents types de parall\'elo\`edres
primitifs} {\bf 136} (1909) 67--181.

\end{thebibliography}
\end{document}